\documentclass[graybox]{svmult}

\usepackage{amssymb, amsmath, latexsym, proof}
\usepackage{xcolor,soul}

\usepackage{type1cm}        
\usepackage{makeidx}         
\usepackage{graphicx}        
\usepackage{multicol}        
\usepackage[bottom]{footmisc}

\usepackage{newtxtext}       %
\usepackage[varvw]{newtxmath}       

\makeindex             

\newcommand{\bl}{\begin{lemma}}
\newcommand{\el}{\end{lemma}}
\newcommand{\bt}{\begin{theorem}}
\newcommand{\et}{\end{theorem}}
\newcommand{\bcor}{\begin{corollary}}
\newcommand{\ecor}{\end{corollary}}
\newcommand{\bp}{\proof{.}}
\newcommand{\ep}{\eop}
\newcommand{\bpr}{\begin{proposition}}
\newcommand{\epr}{\end{proposition}}
\newcommand{\brem}{\begin{remark} }
\newcommand{\erem}{\end{remark}}
\newcommand{\bd}{\begin{definition} }
\newcommand{\ed}{\end{definition}}
\newcommand{\bex}{\begin{example} }
\newcommand{\eex}{\end{example}}
\newcommand{\beq}{\begin{equation} }
\newcommand{\eeq}{\end{equation}}

\newcommand{\bi}{\begin{itemize}
  }
\newcommand{\ei}{\end{itemize}}
\newcommand{\ben}{\begin{enumerate} }
\newcommand{\een}{\end{enumerate} }

\newenvironment{enumr}{

\begin{enumerate}     }{\end{enumerate}

}

\newcommand{\benr}{\begin{enumr}
  }
\newcommand{\eenr}{
\end{enumr}}

\newcommand{\ignore}[1]{}

\newcommand{\al}[1]{\forall #1\:}
\newcommand{\ex}[1]{\exists #1\:}

\newlength{\hilflh}

\newcommand{\naturals}{\mathbb{N}}

\renewcommand{\emptyset}{\varnothing}

\newcommand{\cF}{{\mathcal F}}

\newcommand{\cP}{{\mathcal P}}

\newcommand{\cA}{{\mathcal A}}
\newcommand{\cM}{{\mathcal M}}

\newcommand{\cO}{{\mathcal O}}

\newcommand{\cS}{\mathcal{S}}

\newcommand{\ga}{\alpha}
\newcommand{\gb}{\beta}

\renewcommand{\ge}{\varepsilon}
\newcommand{\gl}{\lambda}

\newcommand{\gw}{\omega}

\renewcommand{\phi}{\varphi}

\newcommand{\Con}{\mathrm{Con}}

\newcommand{\PRA}{\mathsf{PRA}}

\newcommand{\PA}{\mathsf{PA}}

\newcommand{\ZF}{\mathsf{ZF}}

\newcommand{\Lim}{\mathrm{Lim}}

\newcommand{\iffdef}{\stackrel{\text{def}}{\iff}}

\newcommand{\nat}{\naturals}

\renewcommand{\leq}{\leqslant}
\renewcommand{\geq}{\geqslant}

\newcommand{\eop}{$\Box$ \protect\par \addvspace{\topsep}}

\newcommand{\TI}{\mathrm{TI}}

\begin{document}

\title*{Automatic structures and the problem of natural well-orderings}
\author{Lev D. Beklemishev\orcidID{0000-0002-2949-0600} and\\ Fedor N. Pakhomov\thanks{The work of Fedor Pakhomov was funded by the FWO grant G0F8421N.}}
\institute{Lev D. Beklemishev \at Steklov Mathematical Institute, Gubkina str. 8, 119991 Moscow, Russia, \email{bekl@mi-ras.ru}
\and Fedor N. Pakhomov \at Steklov Mathematical Institute, Gubkina str. 8, 119991 Moscow, Russia \at Ghent University,  Krijgslaan 281 S8, 9000 Ghent, Belgium, \email{fedor.pakhomov@ugent.be}}
%
%
\maketitle


\ignore{
\title{Automatic structures and the problem of natural well-orderings}

\author[1]{Lev D. Beklemishev\footnote{The corresponding author,  e-mail: \texttt{bekl@mi-ras.ru}.}}
\author[1,2]{Fedor N. Pakhomov}

\affil[1]{Steklov Mathematical Institute of Russian Academy of Sciences}
\affil[2]{Gent University, Belgium}

\maketitle}

\abstract{
We explore the idea of using automatic and similar kind of presentations of structures to deal with the conceptual problem of natural proof-theoretic ordinal notations. We conclude that this approach still does not meet the goals.
}


\section{Introduction}
This paper is written for the Festschrift volume dedicated to the 75th anniversary of Johann Makowsky. While thinking on the topic that would be appropriate for this volume, we decided to go for one that would link the interests of the authors in proof theory and some topics that play a role in Johann's own work: classes of structures with decidable (MSO) theories and interpretations.
We decided to record our attempts to explore one particular approach to the (in)famous problem in proof theory --- the problem of canonicity of proof-theoretic ordinal notation systems. Even though this problem is well-known, few written accounts and discussions of it exist in the literature.

Arguably, historically the first encounter with this problem occurs in the famous work of Alan Turing ``System of logics based on ordinals''~\cite{Tur39} presenting the content of his PhD dissertation. Substantial discussions are found in the papers by Georg Kreisel, in particular \cite{KraOZ}, and especially Solomon Feferman~\cite{Fefbug} who describes it as one of the three problems that `bug him'.

The problem requires delineating canonical, well-behaved ordinal notation systems from pathological ordinal notation systems, such as a computable well-ordering of order type $\omega$ such that the associated induction scheme implies the consistency of Peano arithmetic (Kreisel’s example). All examples of such pathological well-orderings have some external notions (such as consistency statements) encoded into them; however, it is unclear exactly what is required from an ordinal notation system to forbid such counterexamples and at the same time to be sufficiently general.

The absence of a good mathematical solution of this problem is not only annoying, but it makes the basic question what constitutes an ordinal analysis of a formal theory intuitive rather than fully rigorous. The commonly used expression `to calculate the proof-theoretic ordinal of a theory’ does not really have a definite meaning: it is not clear in which terms the result of this `calculation' needs to be specified. Thus, ordinal analysis, the characterization of proof-theoretic ordinals of theories, is sometimes described as an art~\cite{Rat06,Rat99}. We usually recognize individual ordinal notation systems arising in concrete situations as natural, the simplest one is the system of ordinal notation for $\ge_0$ based on Cantor normal forms. Many other systems are surveyed, e.g., in~\cite{Fefbug,Poh,Rat99}. However, a general mathematical definition of a `natural ordinal notation system' is lacking.

At this point we would like to remark that the term `natural ordinal notation system' may be somewhat misleading. Our understanding of this problem is purely  mathematical: the goal is not to explicate the vague notion of  `naturality' in a philosophical sense of the word, but rather to find a suitably general definition or a framework that would delineate a wide class of ordinal notation systems (well-ordering representations) suitable for proof-theoretic analysis.


Our main goal in this note is to consider this problem from the point of view of the theory of automatic structures. We will motivate this approach and see where it will lead us to.
Our attempts do not really solve the problem, however we believe that there is value in studying failures. In the process we learn new insights, and ultimately such a study may help to point us in the right direction.

\section{Turing}
Arguably the first discussion of the \emph{problem of natural ordinal notations}, as we would call it today, is found in the work of Alan Turing~\cite{Tur39}. There, Turing studied transfinite progressions of theories based on iteration of the process of extending a theory by consistency assertions and by some more general reflection principles. The main goal of this study was to obtain a classification of arithmetical sentences\footnote{Turing considered sentences of logical complexity at most $\Pi_2^0$ that he called `number-theoretic theorems'.} according to the stages of this process. To quote:

\begin{quote} \emph{
We might also expect to obtain an interesting classification of number-theoretic theorems according to ``depth''. A theorem which required an ordinal $\alpha$ to prove it would be deeper than one which could be proved by
the use of an ordinal $\gb$ less than $\ga$. However, this presupposes more than is justified. }
\end{quote}

The last sentence here apparently indicates that Turing realized that his results fall short of the stated goal.
To explain this, we first remark that to construct a Turing progression
$$T_0:=T,\quad T_{\ga+1}=T_\ga+\Con(T_\ga),\quad T_\gl=\bigcup_{\ga<\gl}T_\ga, \text{if $\gl\in\Lim$},$$
one really needs to associate theories $T_\ga$ with ordinal notations (or some kind of constructive representations) rather than with ordinals $\ga$ in the set-theoretic sense. In order to formulate consistency assertions $\Con(S)$ according to the recipe of G\"odel, the axiom set of a theory $S$ must be r.e.\ and represented in the language of arithmetic by a $\Sigma_1$-formula. Thus, the arithmetical formula $\Con(T_\ga)$ must somehow refer to the ordinal $\ga$, and one has to deal with computable and arithmetized ordinal representations rather than with the ordinals themselves. With this understanding, Turing showed how to accurately define such a progression for a given ordinal notation system. This was further elaborated by Feferman in 1962 \cite{Fef62}. Both Turing and Feferman used Kleene's universal ordinal notation system $\cO$ for this purpose.

One of the main results of Turing's paper can be stated (in modern terms) as follows, where $|a|$ denotes the order type of $a\in\cO$.

\begin{theorem}[Turing]
For each true $\Pi^0_1$-sentence $\pi$ there is an ordinal notation $a\in \cO$ such that $|a|=\gw+1$ and $T_{a}$ proves $\pi$.
\end{theorem}

Turing's completeness result is a mixed blessing. The negative side of it is that any true $\Pi_1^0$-sentence can be proved already at stage $\gw+1$ of a suitable Turing progression. These include all sentences of the form $\Con(T_b)$ with $|b|$ much larger than $\gw+1$, which entails that the theories $T_a$ heavily depend on particular ordinal representations rather than on their order types. Thus, the idea of a meaningful classification of sentences according to the progression stages breaks down. Turing put it in a remarkably pessimistic form:

\begin{quote} \emph{
This completeness theorem as usual is of no value. Although it shows, for instance, that it is possible to prove Fermat's last theorem with $\Lambda_P$ (if it is true) yet the truth of the theorem would really be assumed by taking a certain formula as an ordinal formula}\footnote{$\Lambda_P$ is progression based on iteration of consistency, and `ordinal formulas' are his ordinal notations.}.
\end{quote}

Further in the paper Turing suggests a partial way out, a careful selection of specific ordinal notations.

\begin{quote}\emph{
We can still give a certain meaning to the classification into
depths with highly restricted kinds of ordinals. Suppose that we take a particular ordinal logic $\Lambda$ and a particular ordinal formula $\Psi$ representing the ordinal $\ga$ say (preferably a large one), and that we restrict ourselves to ordinal formulae of the form $\mathrm{Inf}(\Psi,a)$.\footnote{These formulas define
initial segments of $\ga$.} We then have a classification into
depths, but the extents of all the logics which we so obtain are contained in the extent of a single logic.}
\end{quote}

Thus, Turing essentially suggests to consider linearly ordered subsets of $\cO$ defined by specific `highly restricted' ordinal notations. It seems likely that he means here particular notations such as, for example, the one for $\ge_0$ based on Cantor normal forms, in other words, ordinal notation systems we would describe as `natural'. However, he admits that by doing this we have to give up the idea of classification of \emph{all} true $\Pi_1^0$ sentences. We would refer to this proposal below as \emph{``limited Turing's program.''}

Turing did not pursue this suggestion any further than that. This kind of approach has much later reappeared in the works of U. Schmerl~\cite{Schm,Sch82}, it was taken up in~\cite{Bek99b} and other papers with many positive results. Turing's idea has been used to classify, for example, $\Pi_n^0$-consequences of specific theories, such as $\PA$ and its predicative extensions. These results were indeed based on highly specific `natural' ordinal notations.

We know many examples of natural ordinal notation systems for fairly large constructive ordinals. Several types of such notation systems are reviewed in Feferman's paper \cite{Fefbug}. However, we still lack a general understanding of what constitutes a `natural' ordinal notation system. Essentially the same problem appears in the study of proof-theoretic ordinals based on Gentzen's approach.

\section{Proof-theoretic ordinals}

Since Gentzen gave his consistency proof of Peano arithmetic by transfinite induction for the ordinal $\ge_0$, much of the work in proof theory has been the exploration of the relationships between formal theories and well-orderings \cite{Kra68,Tak,Poh,Poh09, Rat99}. A proof-theoretic study of a formal theory usually culminates in the calculation of its \emph{proof-theoretic ordinal}, that is, a well-ordering representing a bound on the strength of the system. With this ordinal, the other important characteristics of a formal theory, such as its class of provably total computable functions, are often connected.

Although it is not always duly emphasised in proof-theoretic literature, there are different ways of associating ordinals to theories. They are sensitive to different levels of logical complexity and lead, in general, to inequivalent notions of proof-theoretic ordinals. The most common notion, known since Gentzen and prevalent in the work on proof theory, is the so-called \emph{$\Pi_1^1$-ordinal} which is defined as the supremum of order types of recursive well-orderings that are provably well-founded in a given formal theory. Of course, for this definition to be applicable, the language of the theory must be able to define recursive relations and to express the well-foundedness property of such relations (which is $\Pi_1^1$-complete). Thus, we usually assume the language to contain at least that of first-order arithmetic with free second-order variables (denoting arbitrary sets of natural numbers). A fundamental result of Gentzen stated in these terms is that the $\Pi_1^1$-ordinal of $\PA(X)$, a version of Peano arithmetic in the language expanded by free set variables, is $\ge_0$.

The $\Pi_1^1$-ordinal is a well-defined and robust measure of proof-theoretic strength, however it is not very sensitive: Extension of a theory by true $\Sigma_1^1$-axioms does not change its $\Pi_1^1$-ordinal~\cite{KraOZ,Rat99}. In particular, theories with the same $\Pi_1^1$-ordinal may have vastly different consistency strength and the classes of provably total computable functions. For example, $\PA$ and $\PA+\Con(\ZF)$ have the same $\Pi_1^1$-ordinal. Moreover, Gentzen's consistency proof for  Peano arithmetic is not exactly captured by the characterization of its $\Pi_1^1$-ordinal: Whereas a proof of $\Con(\PA)$ by transfinite induction on $\ge_0$ entails that the well-foundedness of $\ge_0$ is unprovable in $\PA(X)$ (by G\"odel's second incompleteness theorem), the converse cannot in general be concluded.

Attempts to define the notions of proof-theoretic ordinals relevant for
arithmetical complexity classes such as $\Pi_2^0$ and $\Pi_1^0$ only succeed provided some natural system of  ordinal notation is given (as in Turing's limited program),  see~\cite{Bek99b} for a discussion of various proposals of this kind. Without such an assumption, these attempts invariably fail. This issue was treated by Georg Kreisel who constructed a number of pathological counterexamples in order to demonstrate such failures.


\section{Pathological well-orderings} \label{Kreisel}
Given an arithmetically definable binary relation $\prec$ and a class of formulas $\cF$ we denote by  $\mathrm{TI}(\prec;\cF)$ the schema of transfinite induction for $\phi\in\cF$:
$$\al{x}(\al{y\prec x} \phi(y)\to \phi(x))\to \al{x}\phi(x). $$
$\mathrm{TI}(\prec;\cF)$ is true in $\nat$ iff $\prec$ is well-founded w.r.t.\ $\cF$-definable sets. As remarked by Kreisel, Gentzen's consistency proof for $\PA$ was naturally formalizable in Primitive Recursive Arithmetic $\PRA$ extended by $\mathrm{TI}(\prec;\Delta_0)$, where $\prec$ is the canonical primitive recursive well-ordering of order type $\ge_0$, and $\Delta_0$ is the class of primitive recursive arithmetical formulas.

Now we turn to pathological well-orderings.
As our starting point we consider the simplest example from~\cite{Kra68}.

\bpr For any true $\Pi_1^0$
sentence $\pi$, there is a primitive recursive well-ordering $\prec_\pi$ of order type $\gw$ such that $\PRA + \mathrm{TI}(\prec_\pi;\Delta_0)\vdash \pi$.
\epr

\bp\ If $\pi$ has the form $\al{x}\pi_0(x)$ with $\pi_0\in\Delta_0$, then one can define:
$$x \prec_\pi y \iffdef (x < y \land \al{z < x}\pi_0(z)) \lor (y < x \land \ex{z < y}\neg\pi_0(z)).$$
It is easy to see that, once $\pi$ is false and $n$ is the minimal natural number such that $\neg\pi_0(n)$, the ordering $\prec_\pi$ has a $\Delta_0$-definable subset $\{y: y>n\}$ without the least element. Hence $\mathrm{TI}(\prec_\pi;\Delta_0)$ is false. Formalizing this in $\PRA$ shows $\PRA + \mathrm{TI}(\prec_\pi;\Delta_0)\vdash \pi$. On the other hand, since $\pi$ is true in the standard model of arithmetic, the witness of $\neg\pi$ does not exist and $\prec_\pi$ is, in fact, isomorphic to $\gw$ (but not provably so unless one can prove $\pi$).
\ep

This effect is essentially the same as the one observed in the above quotation from Turing: If $\pi$ is Fermat's last theorem (or, better now, some $\Pi_1^0$-equivalent of Riemann's Hypothesis), the truth of $\pi$ is encoded in the fact that $\prec_\pi$ is a well-ordering.
This rather cheap trick is to be compared with Gentzen's theorem stating that $\Con(\PA)$ is provable in  $\PRA + \mathrm{TI}(\prec;\Delta_0)$ where $\prec$ is Cantor's canonical notation system for the ordinal $\ge_0$, but not any proper initial segment of it.

Kreisel's example immediately shows the inadequacy of the naive definition of a proof-theoretic ordinal of a theory $T$ as the order type of the shortest primitive recursive well-ordering $\prec$ such that $\mathrm{TI}(\prec,\Delta_0)$ over $\PRA$ proves $\Con(T)$: Take $\Con(T)$ for $\pi$ and observe that the ordinal of $T$ then equals $\gw$ irrespectively of $T$.


\brem
Kreisel's example shows that $\mathrm{TI}(\prec,\Delta_0)$ can be pathologically strong. However, there are other exmaples showing that this schema can be pathologically weak~\cite{Kra68,Bek00}. For example, Kreisel demonstrated that there are non-wellfounded  primitive recursive linear orderings $\prec$ such that $\PA$ proves the schema of transfinite induction on $\prec$ for arbitrary arithmetical formulas $\mathrm{TI}(\prec,\Pi^0_\infty)$. Accordingly, for all $\ga<\gw_1^{CK}$, one can construct a primitive recursive well-ordering $\prec$ of order type $\ga$, such that $\mathrm{TI}(\prec,\Pi^0_\infty)$ is provabe in $\PA$.  This phenomenon is closely related to the construction of linear orderings well-founded w.r.t.~hyperarithmetical sets, but not actually well-founded, see~\cite{Har68}.
\erem

Other important uses of ordinals in proof theory occur in the constructions of subrecursive hierarchies~\cite{Ros84}. In turn, the latter are used to characterize provably total computable functions of theories~\cite{BW} and as technical tools to prove the independence of certain combinatorial theorems~\cite{KS81}.

An important hierarchy of functions is the so-called \emph{fast-growing hierarchy} $(F_\ga)_{\ga<\Lambda}$ (also known as the extended Grzegorczyk hierarchy) \cite{LW}:
\begin{enumerate}
    \item $F_0(x)=x+1$,
    \item $F_{\alpha+1}(x)=\underbrace{F_\alpha(\ldots (F_\alpha}\limits_{\mbox{\footnotesize $x$-times}}, (x))\ldots)$
    \item $F_{\lambda}(x)=F_{\lambda[x]}(x)$, for limit ordinals $\lambda$.
\end{enumerate}
Here $\Lambda$ is supposed to be a countable ordinal, and $\gl[\cdot]$, for each limit ordinal $\gl<\Lambda$, denotes a \emph{fundamental sequence}, that is, a monotonically increasing sequence of ordinals $\gl[n]$ such that  $\mathsf{sup}_{n<\omega} \lambda[n]=\lambda$. Once a system of fundamental sequences is fixed, the functions $F_\ga$ are uniquely defined for each $\ga<\Lambda$.

If we want the functions $F_\alpha$ to be computable, then $\Lambda$ should be represented as a computable well-ordering; the set $\Lambda\cap\mathsf{Lim}$ and the functions $x+1$ and $x[y]$ should be computable as well.

The typical expectation for non-pathological computable systems of ordinal notations and fundamental sequences is that the growth rate of $F_\alpha$ reflects the size of $\alpha$, e.g., $F_\omega$ is expected to be of the growth rate of the Ackermann function. The fast-growing hierarchy is often used to characterize provably total computable functions of theories. For  natural ordinal notation systems, the functions $F_{\alpha}$ are provably total in $\mathsf{PA}$, if $\alpha<\varepsilon_0$, while $F_{\varepsilon_0}$ grows faster than any $\mathsf{PA}$-provably total computable function.

Examples of computable ordinal notation systems leading to pathological subrecursive hierarchies are well-known and going back at least to~\cite{Fef62a} (see also the references therein). We give an easy example of this kind here.

\bpr For any total computable function $g$ there is a (polytime) computable well-ordering
of order type $\omega+1$ and a  fundamental sequence for $\omega$ such that the growth rate of $F_\omega$ is faster than $g$.
\epr
\bp\ Without loss of generality we switch from $g$ to a perhaps faster growing monotone function $f$ such that $f(x)$ is computable in time polynomial in $f(x)$. The domain of the ordering $\prec$ is the union of $\{\omega\}$ and all  pairs of natural numbers $(n,m)$ such that $m\le f(n)$. We put, for each $n,n',m,m'$, $(n,m)\prec \omega$ and
$$(n,m)\prec (n',m')\text{ if either $n<n'$ or $n=n'$ and $m> m'$.}$$
The fundamental sequence is $\omega[n]:=(n,0)$. Notice that the comparison relation, the domain of the ordering, and the functions $x+1$ and $\omega[x]$ are computable in polynomial time. The only non-trivial algorithm here is the one for the domain.

Since $f(n)$ is computable in time polynomial in $f(n)$, there is a polynomial bound $s(m)$ such that whenever the computation of $f(n)$ does not terminate in $s(m)$ steps,  $f(n)\geq m$. To check if $(n,m)$ belongs to the domain, on input $(n,m)$ run the computation for $f(n)$ for $s(m)$ steps.
If $f(n)$ has not yet terminated tell that $(n,m)$ is in the domain, otherwise check the inequality $m\le f(n)$ directly.

We have $F_\omega(x)\geq f(x)$, since $F_\omega(x)=F_{(x,0)}(x)$ and the position of $(x,0)$ in the ordering is $x+\sum_{y\le x}f(y)\geq f(x)$.
\ep

\section{Complexity considerations} \label{complexity}
How can one define a general class of ordinal representations excluding pathological examples like Kreisel's? We first remark that the idea of restricting the well-orderings to low complexity classes does not really work. The complexity of the formula $\prec_\pi$ is already low ($\Delta_0$, which corresponds to the linear time hierarchy). In fact, one can lower the complexity of $\prec_\pi$ even further.

Let us call a class of binary relations a \emph{basis of r.e.\ sets}, if it is closed under boolean operations and explicit transformations, and any r.e.\ set can be obtained as a projection from a relation in that class. (The term is due to R.~Smullyan, see~\cite{Smu61,Lew79}.) The class of $\Delta_0$-definable relations, for example, is a basis of r.e.\ sets, as is the class of linear time computable relations (on multitape TMs).

Given a basis of r.e.\ sets $\cS$, any $\Pi_1^0$-sentence $\pi$ can be represented in the form $\al{x}\pi_0(x)$ with $\pi_0\in \cS$ in the standard model of arithmetic. We say that $\cS$ is a basis of r.e.\ sets \emph{provably in $T$} if this holds in all models of $T$. If $T$ is sufficiently strong (as strong as $\PRA$ would certainly suffice) we have that linear time computable relations, for example, are a provable basis of r.e.\ sets in $T$.

Let $f$ be a monotone, provably total computable function in $T$ whose graph is linear time computable, and let $g(x):=\min\{y:f(y)>x\}$. Notice that $f$ can be fast growing (e.g.\ as fast as any primitive recursive function if $T=\PRA$), and $g$ is therefore very slow growing. We can modify Kreisel's example as follows:
$$x \prec_\pi y \iffdef (x < y \land \al{z < g(x)}\pi_0(z)) \lor (y < x \land \ex{z < g(y)}\neg\pi_0(z)).$$
We think of $x$ and $y$ now as of binary strings, and $<$ corresponds to the lexicographic ordering. This defines a well-ordering of order type $\gw$ whose computational complexity is only slightly above the linear function of $\max(|x|,|y|)$ assuming that $\pi_0(z)$ is computed in $O(|z|)$ steps. However, the argument for $$\PRA + \mathrm{TI}(\prec_\pi;\Delta_0)\vdash \pi$$ goes through. We remark that one needs the provable totality of $f$ in order to show the existence of a decreasing chain in $\prec_\pi$: If $m$ is the first $z$ such that $\pi_0(z)$ then the decreasing chain starts from $f(m)$.

We conclude that well-ordering representations that do not allow for a version of Kreisel's trick (of incorporating the consistency into the very definition of the ordering) must be more restrictive than most natural complexity classes.

The first class that comes to mind which is not a basis of r.e.\ sets is the class of regular languages. As is well-known, regular languages are closed under boolean operations and projection. So, projection of a regular language cannot be $\Sigma_1$-complete. This leads us to the more general idea of using   automatic presentations of well-orderings as candidates for canonical ones. 

\section{Automatic structures}\label{automatic_structures_section}

Automatic presentations of first-order structures emerged in the fundamental works of B\"uchi, Rabin and others, and are still an ongoing topic of active research. An important milestone in this development was the work by Khoussainov and Nerode~\cite{KhN} where a general program to study such structures was initiated (the notion of automatic structure already appeared in an earlier paper by Hodgson~\cite{Hod82}). A (relational) structure
is called \emph{automatic}  if it has an isomorphic copy where the  universe is a regular set of words and all relations can
be recognized by synchronous multi-tape automata. We recommend the surveys~\cite{KhM07,BG04} for an introduction to the topic of automatic structures and a historical overview.

The study of automata presentable structures was mainly motivated by the problems in computable model theory. The main advantage of automatic structures compared to (polynomially) computable  structures is that they enjoy, in general, nicer computational properties. For example, automatic structures have decidable first-order theories. The class of automatic structures is closed under first-order interpretations (hence under factorization by definable congruences and  definable substructures) and finite products. They also have nice alternative characterizations in terms of logical languages. The following theorem summarizes the results by B\"uchi, Bruy\`ere, Blumensath and Gr\"adel (see \cite{BG04}), where $\cP^{<\gw}(\nat)$ denotes the set of all finite subsets of $\nat$:

\begin{theorem} Let $\cM$ be a first-order structure. $\cM$ is automatic iff any of the following conditions hold:
\benr
\item $\cM$ is first-order interpretable in $(\nat,\cP^{<\gw}(\nat); \in, S)$, where $S$ is the successor function.\footnote{We consider here the two sorted (weak second-order) structure as a first-order structure in the usual way. This allows us to use the standard notion of first-order interpretation here (see e.g.~\cite{BG04,Mak04}).}
\item $\cM$ is first-order interpretable in B\"uchi arithmetic $(\nat;+,V_2)$, where $V_2(x)$ is the function that returns the maximal power of $2$ dividing $x$.
\item $\cM$ is first-order interpretable in $(\{0,1\}^*; \sqsubset,S_0,S_1,E)$, where $\sqsubset$ is the prefix relation, $S_i$ are successor relations, and $E(x,y)$ is true iff the words $x,y$ have equal length.
\eenr
\end{theorem}

\brem
In (ii) above the function $V_2$ can be replaced by any other $V_k$, with $k>2$, as all structures of the form $(\nat;+,V_k)$ are mutually interpretable. This is possible because first order interpretations allow for relativization of quantifiers. Similarly, every automatic structure is automatic over the binary alphabet. In contrast, it is well-known that $V_n$ is definable in $(\nat;+,V_k)$ iff $n$ and $k$ are multiplicatively dependent, that is, if $n^m=k^l$, for some $k,l$~\cite{Vil92}.
\erem

Automatic well-orderings have been studied early on, and one of the basic results is the following theorem by Delhomm\'e~\cite{Del04}:

\bt \label{auto} An ordinal $\ga$ has an automatic presentation iff $\ga<\gw^\gw$. Moreover, from an automaton recognizing the binary relation one can effectively construct a Cantor normal form presentation of its order type $\ga$:
$$\ga=\gw^{n_1}m_1+\gw^{n_2}m_2+\cdots+\gw^{n_k}m_k$$
with $m_i>0$ and $n_1>n_2>\cdots >n_k\geq 0$.
\et

\bcor The isomorphism problem for automatic well-orderings is decidable. \ecor

This result shows that automatic presentations of well-orderings are equivalent to Cantor normal form presentations (of ordinals below $\gw^\gw$). The theorem itself is provable by elementary methods, which means that automatic well-orderings are natural in the sense of proof theory. However, a major drawback is that they do not work beyond the very small ordinal $\gw^\gw$, which is at the lowermost end of ordinal notations of interest in proof theory (it corresponds to the proof-theoretic ordinal of primitive recursive arithmetic).

Can one do better? A way out is to look for more general types of presentations than the automatic ones. In the recent literature several other kinds of automatic-like  structures were considered, especially: \emph{tree-automatic} (which accept as inputs finite labelled trees), \emph{B\"uchi automatic}  (which accept as inputs infinite words), and \emph{Rabin automatic} ones (which accept as inputs infinite binary labelled trees). A structure is called tree (respectively, B\"uchi, Rabin) automatic, if its domain and basic relations are recognizable by tree (respectively, B\"uchi, Rabin) automata. Classical results of  B\"uchi~\cite{Buc} and Rabin~\cite{Rab69} relate this to definability in (weak) monadic second order structures of, respectively, natural numbers with successor and the infinite binary tree. See also \cite{BG04}.

\begin{theorem}[B\"uchi, Rabin]
\benr
\item A structure is B\"uchi automatic iff it is definable in $(\nat,\cP(\nat);\in,S)$.
\item A structure is tree-automatic iff it is definable in \\ $(\{0,1\}^*,\cP^{<\gw}(\{0,1\}^*);\in,S_0,S_1)$.
\item A structure is Rabin automatic iff it is definable in  \\ $(\{0,1\}^*,\cP(\{0,1\}^*);\in,S_0,S_1)$.
\eenr
\end{theorem}

The class of tree-automatic structures is genuinely larger than the class of automatic ones. For example, the structure of natural numbers with multiplication (Skolem arithmetic) is tree-automatic but not automatic. In the case of well-orderings, however, the gains are relatively minor. For tree-automatic ordinals, Delhomm\'e~\cite{Del04} showed that Theorem \ref{auto} can be improved by one exponent: An ordinal $\ga$ is tree-automatic iff $\ga<\gw^{(\gw^\gw)}$.

The classes of B\"uchi and Rabin automatic structures are advantageous when one wants to represent  uncountable structures. Their first-order theories are decidable and in general the classes have nice closure properties. For example, the structure of reals with addition is B\"uchi automatic~\cite{BG04}.

Proof-theoretic ordinal notation systems are supposed to be interpreted in first-order arithmetic (otherwise, transfinite induction schema can hardly be stated), and in this case we want ordinal representations to be finite objects.   However, e.g.\ in the case of the binary tree, if one only considers MSO definable structures whose elements correspond to \emph{finite} sets of words, then these will be the same as tree-automatic (or WMSO definable) structures. Thus, prima facie there is no advantage in considering B\"uchi and Rabin representations of well-orderings. However, this discussion brings us to a few useful general observations.

Firstly, the class of relations definable in decidable theories, such as MSO theories of the binary tree, can never be a basis of r.e.\ sets in the sense of Section \ref{complexity}. Hence, they are good candidates for avoiding pathological well-ordering counterexamples. The more expressive is such a decidable theory, the larger a class of interpretable well-orderings we potentially obtain. This draws our attention towards expressively strong decidable theories.

Secondly, suppose a structure $\cA$ has a decidable monadic second-order theory (as opposed to its first-order theory). If a relation $\prec$ is interpreted in $\cA$ in such a way that its domain consists of elements of $\cA$ (is MSO definable), then the well-foundedness of $\prec$ is expressible by the MSO formula
$$\al{X}\al{x}(\al{y\prec x}y\in X\to x\in X)\to \al{x}x\in X.$$
Hence, given a formula defining $\prec$, one can effectively check whether $\prec$ is a well-ordering. Thus, the class of interpretable well-orderings will have an additional nice property of being effectively recognizable.

Thirdly, because of the previous property, if a model has a decidable MSO theory in the language expanded by constants for each element of the domain, then one can enumerate all the well-orderings definable in it. Hence, their order types will be uniformly bound by some constructive ordinal $\ga<\gw_1^{\text{CK}}$.

The only example we know of where this idea works beyond the tree-automatic ordinals is the so-called Caucal hierarchy.

\section{Caucal hierarchy}
Muchnik's theorem~\cite{Sem84} generalizes the decidability result of Rabin by establishing that the decidability of the MSO theory of a structure is preserved under the iteration operation that maps relational structures $\mathcal{M}=(M;R_0,\ldots,R_{n-1})$ to certain natural structures $\mathcal{M}^*$ on finite words over its domain $M$. In particular, a well studied corollary of Muchnik's results is the decidability of MSO theories of graphs from the Caucal hierarchy~\cite{Cau96,Cau02}.

Here we consider directed graphs with edges colored in finitely many colors such that, for all vertices $x,y$ and color $i$, there is at most one edge from $x$ to $y$ of that color.  Formally, we say that a structure $\mathcal{G}$ is a \emph{directed graph with colored edges} if its signature consists of finitely many binary relations $\{R_i\mid i< n\}$. Later in this section we will simply call them graphs.

We have two natural operations on these graphs preserving the decidability of MSO theories:
\begin{enumerate}
    \item MSO interpretations\footnote{Transductions in the terminology of \cite{Mak04}.}, i.e., the interpretations where the first-order domain is interpreted by a one-dimensional MSO definable set of first-order elements and the relations are interpreted as MSO definable relations.
    \item The unfolding operation $\mathsf{Unf}$, where for a graph $\mathcal{G}$, the vertices of $\mathsf{Unf}(\mathcal{G})$ are all possible paths $\mathsf{Path}(\mathcal{G})$ through $\mathcal{G}$, and we have an $i$-colored edge from a path $\alpha=(v_0,\ldots,v_m)$ to $\beta=(u_0,\ldots,u_k)$ if $k=m+1$, $v_0=u_0$,$\ldots$,$v_m=u_m$, and in $\mathcal{G}$  there is an $i$-colored edge from $v_m$ to $u_{m+1}$.
\end{enumerate}
Level $0$ of the Caucal hierarchy $C_0$ consists of all finite graphs, and level $n+1$ of the hierarchy $C_{n+1}$ consists of all the graphs that are MSO  interpretable in the unfoldings of graphs from $C_n$.

From the automata-theoretic perspective, the Caucal hierarchy naturally corresponds to \emph{higher-order pushdown automata} first introduced by Maslov~\cite{Mas74}. We call a \emph{$0$-pds} over a finite alphabet $A$ (pds stands for \emph{higher order pushdown store}) just a letter from $A$. A  \emph{$(n+1)$-pds} over $A$ is a finite sequence of $n$-pds. We think about $(n+1)$-pds as a stack of $n$-pds. Further, we have a natural $\mathsf{pop}^k$ operation, $0<k\le n$, removing the topmost $(k-1)$-pds (if $k<n$ consider the topmost $(n-1)$-pds and apply $\mathsf{pop}^k$ to it).  The operation $\mathsf{push}^k(a)$ modifies the topmost $k$-pds by creating a copy of the topmost $(k-1)$-pds and putting it on the top of the $k$-pds, followed by replacing the topmost $0$-pds with the letter $a$.

Now one can naturally define the notion of a \emph{$n$-pushdown automaton} as a finite state transition device that uses a single $n$-pds as its memory. We do not go into the details here, but it is important to allow $\varepsilon$-transitions, i.e., the transitions where the automaton can  perform an operation on $n$-pds and transit to a new state,  but is not reading any characters from the input. The \emph{configuration graph} of a (possibly non-deterministic) $n$-pushdown automaton is the graph whose vertices are the states reachable from the starting state (with an empty $n$-pds). Its edges are labeled with $\varepsilon$ and characters of the input alphabet and correspond to the one-step transitions. As proved by Carayol and W\"ohrle~\cite{CW03}, level $n$ of the Caucal hierarchy consists precisely of (the graphs isomorphic to) $\varepsilon$-contractions of the configuration graphs of $n$-pushdown automata.

In the previous section we connected ordinals with structures with decidable MSO theories by means of first-order interpretations. Although this question is also natural to ask for the case of the Caucal hierarchy, to the best of our knowledge it is open.

In this section we will be looking at a stronger notion of representability of ordinals in the sense of being MSO interpretable. Structures in the Caucal hierarchy have decidable MSO theories and all their  elements are definable. Hence, the remarks at the end of the previous section apply. In particular, one can effectively recognize well-orderings in the Caucal hierarchy, and there is a uniform bound on their order types.
The following theorem by Braud and Carayol~\cite{Bra09,BraCar10} explicitly describes this bound.

Define $\gw_0:=1$ and, for $k<\gw$, let $\gw_{k+1}:=\gw^{\gw_k}$. Then $\ge_0$ is the supremum of $\{\gw_k:k<\gw\}$.

\begin{theorem} An ordinal $\alpha$ is MSO interpretable in a graph from the Caucal hierarchy iff $\alpha <\varepsilon_0$. An ordinal $\ga$ is interpretable in the $n$-th level of the Caucal hierarchy iff $\alpha<\gw_{n+1}$.
\end{theorem}

\brem We do not really know if the Cantor normal form of an ordinal can be computed from its Caucal hierarchy representation. Nor do we know if the isomorphism problem for well-orderings in the Caucal hierarchy is decidable.
\erem

\brem
We know that MSO theories of the structures in the Caucal hierarchy are decidable. A fortiori, this holds for the well-orderings in the Caucal hierarchy.

By the results of L. Ko{\l}odziejczyk and H. Michalewski~\cite{KM16}, Rabin's theorem on the decidability of the MSO theory of the binary tree is surprisingly strong: it is unprovable in the fragment of the second-order arithmetic with $\Delta^1_3$-comprehension axioms. This theory is very much stronger than $\PA$ and its proof-theoretic ordinal is larger than all currently known ones. We do not really know the corresponding lower bound for the Caucal hierarchy, but in any case it can only be worse.

Therefore, we are in a curious situation that Caucal hierarchy presentations, in general, may not be recognizably decidable within a given formal system. It puts some doubts on the hope that any Caucal presentation of a well-ordering would be provably isomorphic to a computable one (such as the canonical one based on Cantor normal forms). However, there is a caveat: By a well-known theorem of B\"uchi~\cite{Buc65,Buc73} any countable well-order has a decidable MSO theory. This theorem seems to be proof-theoretically weaker than Rabin's.
So, Caucal presentations of \emph{well-orderings} may not actually require as strong axioms as $\Delta^1_3$-comprehension.

In any case, there is a strange discrepancy between the order types of Caucal representable well-orderings and the very strong axioms needed to show that they are decidable.
\erem

\section{Fundamental sequences in the Caucal hierarchy}
As we discussed in Section \ref{Kreisel}, the assignment of fundamental sequences to ordinals within an ordinal notation system is important for the construction of subrecursive hierarchies of functions such as the fast-growing hierarchy. Given that such hierarchies play a significant role in proof theory, in this section we will discuss a result from \cite{Pakh15b} stating that fundamental sequences can also be represented with the Caucal hierarchy and yield, under some natural conditions, the hierarchies of functions equivalent to the one for the standard fundamental sequences assignment. We are very sketchy and refer the reader to \cite{Pakh15b} for more details.

Firstly, we remark that an ordinal represented in the Caucal hierarchy can always be expanded by a system of fundamental sequences in the same hierarchy. We represent a system of fundamental sequences by the predicate $\mathsf{FS}(\beta,\gamma)$ expressing that $\beta$ is an element of the fundamental sequence for $\gamma$. Here we naturally assume that only the limit ordinals have non-empty fundamental sequences and, for any limit ordinal $\lambda$, the set of ordinals $\{\beta\mid \mathsf{FS}(\beta,\lambda)\}$ has order type $\omega$ and $\lambda$ as its limit. We then let $\lambda[n]_{\mathsf{FS}}$ denote the $n$-th element of the set $\{\beta\mid \mathsf{FS}(\beta,\lambda)\}$.

To see that, given a well-ordering in the Caucal hierarchy, some system $\mathsf{FS}$ exists that can be defined in the Caucal hierarchy, we consider deterministic trees. A \emph{deterministic tree} is a tree where from each vertex there is at most one outgoing edge of any given color. It is known~\cite{CW03} that every graph from $C_n$ is MSO interpretable in a deterministic tree from $C_n$. It is then fairly easy to expand an MSO interpretation of an ordinal in a deterministic tree by an MSO definable system of fundamental sequences $\mathsf{FS}$.

A commonly considered nice property of systems of fundamental sequences is the so-called \emph{Bachmann property}~\cite{Bach59,Schm77,Ros84}. It demands that each $\lambda[\cdot]\colon \omega\to\lambda$ is strictly increasing and that for any limit $\alpha$ from $(\lambda[n],\lambda[n+1]]$ the value $\alpha[0]\geq \lambda[n]$. As explained in \cite{Pakh15b}, if a well-ordering equipped with a system of fundamental sequences $\mathsf{FS}$ is MSO interpreted in a deterministic tree, then one can construct another MSO definable system $\mathsf{FS}'\subseteq \mathsf{FS}$ with the Bachmann property.

As before, when we have a computable ordinal equipped with system of computable fundamental sequences $\cdot[\cdot]_{\mathsf{FS}}$, we can form the associated fast-growing hierarchy.
If $\mathsf{FS}$ satisfies the Bachmann property, then the asymptotic growth rate of the functions $F_\alpha$ from the hierarchy grows monotonically in $\alpha$~\cite{Ros84}. One of the criteria of an ordinal notation system of being natural is that it leads to the levels of fast-growing hierarchy having the expected growth rate reflecting the value of ordinal. The following theorem~\cite{Pakh15b} confirms that this is the case for the Caucal hierarchy representations.

\begin{theorem} \label{Caucal_F} Suppose a structure $(\Lambda,<,\mathsf{FS}_1,\mathsf{FS}_2)$ in the Caucal hierarchy is an ordinal with two systems of fundamental sequences, both satisfying the Bachmann property. For each $i\in\{1,2\}$, let $(F^{\mathsf{FS}_i}_\alpha\colon \mathbb{N}\to\mathbb{N})_{\alpha<\Lambda}$ be the fast-growing hierarchy according to $\mathsf{FS_i}$. Then for each $\alpha<\beta<\Lambda$ there is an $N$ so large that
$$F^{\mathsf{FS}_1}_\alpha(x)<F^{\mathsf{FS}_2}_\beta(x)\text{, for $x\geq  N$.}$$
\end{theorem}
Hence, all systems of fundamental sequences in the Caucal hierarchy yield fast-growing hierarchies of similar growth rates.  An important tool to get bounds on the hierarchy functions is the pumping lemma for higher-order pushdown automata proved by Parys~\cite{Par12}.

The standard system of fundamental sequences for ordinals $\gl<\ge_0$ given in Cantor normal form is defined by:
$$
\gl[n] := \begin{cases}
 \ga+\gw^\gb\cdot (n+1), &\text{ if $\gl=\ga+\gw^{\gb+1}$}; \\
\ga+\gw^{\gb[n]}, & \text{ if $\gl=\ga+\gw^{\gb}$ and $\gb\in\Lim$}.
\end{cases}
$$
The standard system satisfies the Bachmann property and is MSO definable in $(\Lambda,<)$ for ordinals $\Lambda<\omega^\omega$ (but not above). Hence, we obtain the following corollary.

\bcor \label{Caucal_F_std} Suppose $(\Lambda,<,\mathsf{FS})$ in the Caucal hierarchy is an ordinal $<\omega^\omega$ together with a system of fundamental sequences satisfying the Bachmann property. Let $(F^{\mathsf{FS}}_\alpha\colon \mathbb{N}\to\mathbb{N})_{\alpha<\Lambda}$ be the fast-growing hierarchy according to $\cdot[\cdot]_{\mathsf{FS}}$ and let $(F_\alpha\colon\mathbb{N}\to\mathbb{N})_{\alpha<\omega^\omega}$ be the fast-growing hierarchy (defined using the standard fundamental sequences).   Then for each $\alpha<\beta<\Lambda$ there is an $N$ so large that
$$F_\alpha(x)<F^{\mathsf{FS}}_\beta(x)\text{\;\;  and  \;\;}F^{\mathsf{FS}}_\alpha(x)<F_\beta(x)\text{, for $x\geq N$.}$$
\ecor

\ignore{
Here we give a sketch of a proof that uses ideas broadly similar to \cite{Pakh15b}, but not exactly the same.

The finite levels of two hierarchies of course coincide with the standard fast-growing hierarchy, since there the definition doesn't rely on fundamental sequences at all.

To study levels of the hierarchy $\geq\omega$ we consider an interpretation of $(\Lambda,<,\mathsf{FS}_1,\mathsf{FS}_2)$  in a deterministic tree $\mathcal{T}=(T,R_0,\ldots,R_{m-1})$ from $n$-th level of Caucal hierarchy. That is the ordinals $<\Lambda$ are represented by the words in the alphabet $\{0,\ldots,m-1\}$ (vertices in $\mathcal{T}$ correspond to words given the path from the root of $\mathcal{T}$ to respective vertices). The desired inequalities about the hierarchies are easy to prove using the primitive recursiveness of  the functions $\alpha-1$ (on successor ordinals $\alpha<\Lambda$), $\lambda[x]_{\mathsf{FS}_i}$ (on limit ordinals $\lambda<\Lambda$ and naturals $x$) and $f\colon \Lambda\cap\mathsf{Lim}\times\mathbb{N}\to\mathbb{N}$
$$f(\lambda,k)=\min \{l\mid \lambda[l]_{\mathsf{FS}_2}>\lambda[r]_{\mathsf{FS}_1}\text{, for $r\le k$}\}.$$
Note that any primitive recursive function is bounded from above by $F_k$, for some natural $k$. The general idea of how to get primitive recursive bounds on the functions is to use the pumping lemma for higher-order pushdown graphs proved by Parys  \cite{Par12} (in fact it even gives elementary recursive bounds).
\ep }

We do not know if the counterpart of Corollary~\ref{Caucal_F_std} holds without the restriction on $\Lambda$. One can also ask similar questions about the other classes of automatically represented well-orders considered in Section~\ref{automatic_structures_section}.

\section{Automatic well-founded relations}
We have seen that the notion of automatic well-order provides natural ordinal representations, however is too restrictive for proof-theoretic applications. The situation is the same for well-founded partial orders~\cite{KhM09}. In proof theory, transfinite induction can be stated more generally for well-founded relations rather than just for well-founded orders (or partial orders).

Are automatic presentations of well-founded relations always natural? We answer this question negatively by employing a construction from~\cite{KhM09} who showed that there exist automatic well-founded relations of arbitrary large ordinal rank $<\gw_1^{CK}$.

\bt For each true $\Pi_1^0$-sentence $\pi$ there exists an automatic well-founded binary relation $R_\pi$ such that $\PRA+\TI(R_\pi,\Delta_0)\vdash \pi$.
\et

\bp\ The proof is essentially an application of a construction from \cite[Theorem 1.2]{KhM09} to Kreisel's example. We only sketch it here and refer to \cite{KhM09} for additional details.

The proof relies on several general  facts. Firstly, the configuration graph of a Turing machine is automatic. Vertices of the graph are tuples of words (representing the content and the position of the head on each tape). Edges represent one-step transitions of the machine between the configurations.

Secondly, for each Turing machine there is a reversible three tape Turing machine accepting the same language. For such a machine both the in-degree and the out-degree of the configuration graph are at most one, in fact, the graph is a union of chains that are either finite or of the type of natural numbers. Therefore, the graph is a well-founded relation (of rank $\leq\gw$).

Given a $\Pi_1$-sentence $\pi$, let $\prec_\pi$ denote Kreisel's ordering on $\{0,1\}^*$ (as defined in Section \ref{Kreisel}), and let $\cM$ be the Turing machine computing $\prec_\pi$ in the sense that on input $(x,y)$ it outputs `yes' or `no' depending on whether $x\prec_\pi y$ holds. Kreisel's definition of $\prec_\pi$ can be read as an algorithm for computing $\prec_\pi$. Moreover, the associated Turing machine $\cM$ can be represented in arithmetic in such a way that the result of its computation provably in $\PRA$ meets its specification: $\cM(x,y)=\text{`yes'}$ iff $x\prec_\pi y$. Moreover, we can assume $\cM$ to be ($\PRA$ provably) reversible.

Let $(D,E)$ denote the configuration graph of $\cM$. We define the domain of the automatic structure $A$ as $\{0,1\}^*\cup D$ (the union is disjoint). The relation $R_\pi$ is defined as the union of $E$ (on $D$) and:
\begin{itemize}
\item All pairs $(x,z)$ such that $x\in\{0,1\}^*$ and $z$ is the initial configuration of $\cM$ on input $(x,y)$, for some $y\in\{0,1\}^*$;
\item All pairs $(z,y)$ where $y\in\{0,1\}^*$ and $z$ is the final configuration of an accepting computation of $\cM$ on input $(x,y)$, for some $x\in\{0,1\}^*$.
\end{itemize}
Since $(D,E)$ is automatic, it is easy to see that so is $R_\pi$.

Let $R^*$ denote the transitive closure of $R$.
\bl\ \label{emb}
For all $x,y\in\{0,1\}^*$, if  $x\prec_\pi y$ then $x R_\pi^* y$.
\el

\bp\ Suppose $x\prec_\pi y$, then $\cM(x,y)=\text{`yes'}$. Let $I(x,y)$ denote the initial configuration of $\cM$ on input $(x,y)$, and $F(x,y)$ its final configuration. Then we have a path in $(A,R_\pi)$ from $x$ to $y$: $x R_\pi I(x,y) R_\pi^* F(x,y) R_\pi y$. \ep

Using Kreisel's argument we see that, if $\pi$ is false, then there is a descending sequence w.r.t.\ $\prec_\pi$. By Lemma \ref{emb} it generates a descending sequence w.r.t.\ $R_\pi$, hence $R_\pi$ is not well-founded. Formalizing this in $\PRA$ yields that $\PRA+\TI(R_\pi,\Delta_0)\vdash \pi$.

Now we show that $R_\pi$ is actually well-founded. Since $\pi$ is true in the standard model, $\prec_\pi$ is isomorphic to $\gw$.
Consider a nonempty set $X\subseteq A$. If $X\cap \{0,1\}^*\neq \emptyset$ we first consider the $\prec_\pi$ minimal element $a\in X\cap \{0,1\}^*$. By the construction of $R_\pi$ (second item), the elements $u$ such that $uR_\pi a$ must be final configurations of accepting computations of $\cM(x,a)$, for some $x$. Select any such $x$. Then $x\prec_\pi a$ and, by the minimality of $a$, $x\notin X$. The computation chain of $\cM(x,a)$ is finite. We claim that its minimal element in $X$ (if exists) will be $R_\pi$-minimal, otherwise $a$ will be $R_\pi$-minimal. This is clear for all but the first element of the computation chain.
However, if $v$ is the initial configuration then its only $R_\pi$-predecessor is $x$, but $x$ is not in $X$. So, $v$, if in $X$, must be $R_\pi$-minimal.

If $X\cap \{0,1\}^*= \emptyset$, then $X$ is a non-empty subset of the computation graph $D$, which is well-founded by the reversibility condition. The minimal element of $X$ in $D$ will also be $R_\pi$ minimal in $A$. \ep\

This example shows that automatic presentation of a structure is not always nice. Whether or not crucially depends on the choice of the relations of the structure one assumes to be automatic. So, this brings us back  to the question, what kind of structures are proof-theoretic ordinals?

\section{Conclusion and open questions}

We would like to mention the following questions resulting from our analysis. We think they are interesting irrespectively of the problem of naturality of ordinal notation systems.

-- Are there mathematically interesting classes of structures expressively stronger than the Caucal hierarchy for which the MSO theory (FO theory) is decidable? We would like to find such structures (or to show that they do not exist) for important proof-theoretic ordinals such as the Feferman--Sch\"utte ordinal, the Howard ordinal, etc.

-- Can one relate in an intrinsic way the Caucal hierarchy representations of ordinals below $\ge_0$ and Peano arithmetic, that is, to use them directly for a proof-theoretic analysis of $\PA$? For example, it seems natural to represent in the Caucal hierarchy the infinitary derivation trees arising from $\PA$-proofs and possibly the non-well-founded proofs in a cyclic version of $\PA$~\cite{Sim17}.

-- As we have discussed above, automatic well-founded partial orders have very small ranks, whereas automatic well-founded binary relations do not, but neither do they exclude pathological counterexamples.
Are there natural types of well-founded structures, whose automatic presentations are tame, yet much larger proof-theoretic ordinals are presentable? This is related to the more traditional view of proof-theoretic ordinals as orders equipped with an additional structure. In the proof-theoretic literature there are various proposals as to possible structures (see \cite{KraOZ,Rat99,Fefbug,Tak, Poh09, Bek04}).

\section{Acknowledgements}
The first author would like to thank Andrei Muchnik (1958--2007) with whom he has had several memorable  conversations about the problem of natural proof-theoretic ordinal notations. Curiously, they were not related to MSO definability and Muchnik's own remarkable work playing a role in the present paper. We also thank Iskander Kalimullin and Alexander Shen for helpful remarks and references.

\end{document}